\PassOptionsToPackage{prologue,dvipsnames}{xcolor}

\documentclass[11pt]{amsart}

\usepackage[latin1]{inputenc}
\usepackage[active]{srcltx}
\usepackage{fullpage}
\usepackage{amsfonts,enumerate}
\usepackage{comment, enumerate}
\usepackage[mathscr]{euscript}
\usepackage{epsfig}
\usepackage{latexsym}
\usepackage[all]{xy}
\usepackage{amssymb}
\usepackage{amsthm}
\usepackage{amsmath}
\usepackage{amsxtra}
\usepackage{float}
\usepackage{tikz}
\usepackage{tikz-cd}
\usepackage{amssymb}
\usepackage{amsmath}
\usepackage{amsfonts}
\usepackage{amsthm}
\usepackage{multicol}
\usepackage{mathrsfs}
\usepackage{graphicx}
\usepackage{enumitem}
\usepackage{url}
\usepackage{tabularx}
\usepackage{float}
\usepackage{psfrag}
\usepackage{tikz}
\usetikzlibrary{positioning}
\usepackage{pbox}
\usepackage{tkz-graph}
\usetikzlibrary{snakes,fit,positioning,calc}
\usepackage{bm}

\usepackage[dvipsnames]{xcolor}
\usepackage[tickmarkheight=.5em,textwidth=\marginparwidth,textsize=small]{todonotes}
\usepackage[backref=page]{hyperref}
\hypersetup{citecolor=MidnightBlue,
	colorlinks=true,
	linkcolor=MidnightBlue,
	urlcolor=MidnightBlue}

.tfm
.tfm

\theoremstyle{plain}
\newtheorem{prop}{Proposition}[section]

\newtheorem{thm}[prop]{Theorem}

\newtheorem{lemma}[prop]{Lemma}
\newtheorem{conj}{Conjecture}

\newtheorem{thmA}{Theorem}
\newtheorem*{ProblemA}{Question}

\newtheorem{thmB}{Theorem}

\theoremstyle{defi}
\newtheorem{defi}[prop]{Definition}

\theoremstyle{remark}

\newtheorem{remark}[prop]{Remark}

\numberwithin{table}{section}

\DeclareMathOperator{\End}{End}

\DeclareMathOperator{\Gal}{Gal}

\newcommand{\bfM}{\mathcal M}

\newcommand{\F}{\mathbb F}
\newcommand{\Gm}{\mathbb G}
\newcommand{\bfC}{\mathcal C}

\newcommand{\Om}{{\mathcal{O}}}

\newcommand{\GL}{{\rm GL}}

\newcommand{\HGM}{\mathcal H}

\def\z-spec{t}
\newcommand{\mubb}{\mbox{$\raisebox{-0.59ex}
  {$l$}\hspace{-0.18em}\mu\hspace{-0.88em}\raisebox{-0.98ex}{\scalebox{2}
  {$\color{white}.$}}\hspace{-0.416em}\raisebox{+0.88ex}
  {$\color{white}.$}\hspace{0.46em}$}{}}

\newcommand{\n}{\mathfrak{n}}

\makeatletter

\newcommand{\verbatimfont}[1]{%
	\renewcommand{\verbatim@font}{\ttfamily #1}%
}
\renewcommand{\verbatim@font}{\ttfamily\small}
\makeatother

\def\RR{\mathbb R}

\def\QQ{\mathbb Q}

\def\PP{\mathbb P}

\def\vK{v_{\scalebox{0.6}{$\kern-0.04em K$}}}

\DeclareFontFamily{U}{mathb}{\hyphenchar\font45}
\DeclareFontShape{U}{mathb}{m}{n}{
	<5> <6> <7> <8> <9> <10> gen * mathb
	<10.95> mathb10 <12> <14.4> <17.28> <20.74> <24.88> mathb12
}{}
\DeclareSymbolFont{mathb}{U}{mathb}{m}{n}
\DeclareMathSymbol{\lefttorightarrow}{3}{mathb}{"FC}

\def\acts{\mathop{\raise 0.13ex\hbox{\protect\scalebox{0.8}[0.9]{$\lefttorightarrow$}}}}

\tikzset{ 
	tck/.style={decoration={ticks,segment length=1pt,amplitude=2pt,raise=-2pt},black!40,decorate}
}

\def\<#1>{{\left\langle{#1}\right\rangle}}

\def\Z{{\mathbb Z}}             
\def\Q{{\mathbb Q}}             

\def\id#1{{\mathfrak{#1}}}      

\newcommand{\vdashv}{\mathrel{\vdash\mkern-4mu\dashv}}

\begin{document}
	
\title{Trivial solutions in generalized Fermat equations}

	\author{Pedro-Jos\'e Cazorla Garc\'ia}
	
	\address{Departamento de Matem\'atica Aplicada, ICAI, Universidad Pontificia Comillas, Madrid, 28015, Spain}  \email{pjcazorla@comillas.edu}
	\thanks{PJCG is supported by the project ``Hipergeometría'' from the Universidad Pontificia Comillas.}

	\author{Ariel Pacetti}
	\address{Center for Research and Development in Mathematics and Applications (CIDMA),
		Department of Mathematics, University of Aveiro, 3810-193 Aveiro, Portugal}
	\email{apacetti@ua.pt}
	\thanks{AP is supported by CIDMA under the Portuguese
		Foundation for Science and Technology (FCT,
		https://ror.org/00snfqn58) Multi-Annual Financing Program for R\&D
		Units, grants UID/4106/2025 and UID/PRR/4106/2025.}

	\author{Lucas Villagra Torcomian}
	\address{Center for Research and Development in Mathematics and Applications (CIDMA),
		Department of Mathematics, University of Aveiro, 3810-193 Aveiro, Portugal}
	\email{lvillagra@ua.pt}
	
	\thanks{LVT is supported by the European Union's Horizon Europe research and innovation programme under the Marie Sk\l odowska-Curie grant agreement No.\ 101270191.}
	
\keywords{}
\subjclass[2020]{11D41, 11F80, 14H25, 14M25}

\begin{abstract} The modular method is a powerful tool to study 
  generalized Fermat equations. Trivial solutions often become an
  obstruction to the method. While using hypergeometric motives to
  construct Galois representations, some of these solutions give rise to newforms with complex multiplication. In
  the present article, assuming a standard conjecture in Darmon's program, we
  prove that asymptotically any newform with complex multiplication cannot be
  related to a non-trivial primitive solution. This allows to discard all such newforms, including those associated to trivial solutions.
  
\end{abstract}

	\maketitle	
		\section{Introduction}

	Since the proof of Fermat's Last Theorem by  Wiles \cite{Wiles}, an increasing amount of interest has been devoted to the study of \emph{generalized Fermat equations}. These are equations of the form 	
	\begin{equation}\label{eq:GFE}
		x^q+y^p+z^r=0,
	\end{equation}
	where  the triple
        $(q,p,r)\in\Z^3_{>0}$ is known as the \textit{signature}. A
        solution $(a,b,c) \in \Z^3$ of \eqref{eq:GFE} is called
        \textit{primitive} if $\gcd(a,b,c)=1$ and \textit{trivial} if
        $abc=0$. From \cite[Theorem 2]{DarmonGranville},
        (\ref{eq:GFE}) has only finitely many non-trivial
        primitive solutions for a fixed signature $(q,p,r)$ satisfying $\frac{1}{q} + \frac{1}{p} + \frac{1}{r}<1$.
	However, Beal's conjecture states the non-existence of non-trivial solutions when $q,p,r\ge3$ are free variables.
	
The approach used by Wiles to prove Fermat's Last Theorem is now 
commonly known as the \textit{modular method}, and is well-suited for
studying solutions of~(\ref{eq:GFE}) when the exponents move along a
line. Restricting to prime exponents, the signatures belong to one of the following four families (up to relabeling of
the variables): $(p,p,p)$\footnote{This case corresponds to Fermat's Last Theorem.}, $(p,p,r)$, $(r,p,r)$ and $(q,r,p)$, where
$q,r$ are fixed and different, and $p$ varies. Of these, the case $(q,r,p)$ is the least understood in the literature.

	Let us briefly discuss how the modular method works. For this, let $\zeta_k$ be a primitive $k$-th root of unity and let $K=\Q(\zeta_q)^+\cdot\Q(\zeta_r)^+$. Finally, let $\id{p}$ be a prime of $K$ above $p$. The modular method proceeds by contradiction and its main steps are the following:
	
	\begin{description}
\item[Step I] To a hypothetical solution $(a,b,c)$ of (\ref{eq:GFE}),
  attach a geometric Galois
  representation \[\rho_{\id{p}}:\Gal(\overline{\Q}/K)\to\GL_2(K_\id{p}),\]
  such that the conductor $\n$ of the residual representation
  $\overline{\rho}_\id{p}$ is independent of the solution $(a,b,c)$.
		
\item[Step II] Prove that $\overline{\rho}_{\id{p}}$ is modular, that is, that there exists a Hilbert modular form $g$ such that $\overline{\rho}_{\id{p}} \simeq \overline{\rho}_{g,\id{p}}$ .
\item[Step III] Prove that $\overline{\rho}_{\id{p}}$ is absolutely
  irreducible. Then, by level-lowering results, we conclude that there exists a Hilbert newform $f$ of parallel weight $2$ and level $\n$ such that   $\overline{\rho}_{\id{p}}\simeq\overline{\rho}_{f,\id{p}}$.
\item[Step IV] Compute $S_2(\n)$ (which, by Step I, is independent of the solution $(a,b,c)$) and show that
  $\overline{\rho}_\id{p}\not\simeq\overline{\rho}_{f,\id{p}}$ for
  all newforms $f$ in $S_2(\n)$.
\end{description}
	
Historically, Step I was performed by associating to a solution
$(a,b,c)$ a rational elliptic curve, so that $\rho_{\id{p}}$ is its
$\id{p}$-adic Galois representation. In more recent articles, following Darmon's pioneering paper \cite{Darmon}, the
elliptic curve was replaced by  a $\GL_2$-type  Jacobian of  a hyperelliptic curve. These curves are referred to as Frey hyperelliptic
curves and have been used to study the families $(p,p,r)$ \cite{Darmon,MR3690598,ChenKoutsianas},
$(r,r,p)$ \cite{BCDFmultifrey,BCDF1,BCDF2}, and $(2,r,p)$ \cite{PedroAngelosLucas}.

A further generalization is due to Golfieri and  the second author  \cite{GP}, where the general signature case was addressed by using hypergeometric motives (as implicitly used by Darmon in
\cite{Darmon}), solving Step I completely. More precisely, the
following holds.

\begin{thm}\label{thm:rep-existence}  Let $q>r>2$  be  prime numbers and $K=\Q(\zeta_q)^+\cdot\Q(\zeta_r)^+$. Let $p>3$ be a prime number and $\id{p}$ a prime ideal of $K$ above $p$. For every hypothetical primitive solution $(a,b,c) \in \Z^3$ to
	\[
	x^q+y^p+z^r=0,
	\]
	there is a Galois representation 
	\[
	{\rho}_\id{p} : \Gal(\overline{\Q}/K) \to \GL_2(K_\id{p}),
	\]
	attached to $(a,b,c)$  satisfying that its residual representation $\overline{\rho}_{\id{p}}$ is unramified outside primes dividing $pqr$ and  finite at $p$.
\end{thm}

The representation $\rho_{\id{p}}$ comes from the hypergeometric motive
$\HGM((\frac{1}{r},-\frac{1}{r}),(\frac{1}{q},-\frac{1}{q})|t)$ (see
\S \ref{sec:set-up-motives} for details).
Regarding Step II, in \cite[$\S 5$ and $\S 6$]{GP} some modularity
results were obtained. For example, the case $r=3$ and $q\geq 5$ is completely solved in Theorem 6.2 in loc.\ citt., while the general case is proved in \cite{Modularity-HGM}.

Step III is particularly difficult to approach in full generality.  In
\cite[Corollary 6.2]{PacettiVillagra}, the irreducibility of
$\overline{\rho}_{\id{p}}$ is proved only under certain divisibility conditions on the solution $(a,b,c)$. 
In order to fulfill Step III, we  assume the veracity of the following ``large-image'' conjecture (Conjecture 4.1 in \cite{Darmon}).
	
\begin{conj}\label{conj:large-image}
  Let $K$ be a totally real field and $L$ any  number field. Then there exists
  a constant $C(L,K)$, depending only on $L$ and $K$, such that for
  any abelian variety $A/L$ of $\mathrm{GL}_2$-type with
  \[
    \End_L(A)\otimes \Q = \End_{\overline{L}}(A)\otimes \Q \simeq K,
  \]
  and all  prime ideals $\mathfrak{p}$ of $K$ with norm  
  larger than $C(L,K)$, the image of the mod $\id{p}$-representation associated to $A$ contains $\mathrm{SL}_2(\F_{\mathfrak{p}})$.
\end{conj}
Assuming Conjecture \ref{conj:large-image}, if the hypergeometric
motive does not have complex multiplication, it follows that
$\overline{\rho}_{\id{p}}$ is absolutely irreducible when the norm of $\id{p}$
is sufficiently large. Then \cite[Theorem 3.2.2]{MR3294620} proves the
existence of the newform $f$ given in Step III (see also
\cite{Jarvis,Rajaei,Fujiwara}).
	
With the first three steps completed, it remains to perform Step IV. For this, there are two different difficulties. The first one is computational. Indeed, when the norm of the conductor $\n$ is large, the dimension of the space of newforms in $S_2(\n)$ makes the computation unfeasible. For instance, this already occurs for $q=7$ and $r=3$.

The second challenge is more theoretical and is related to the
existence of newforms with complex multiplication in the space
$S_2(\n)$ (these forms tend to consistently pass the elimination
process due to the large number of zeroes in their Fourier
coefficients). The \emph{trivial solutions} 
\[
  {(\pm 1, \mp 1,0), \ (\pm 1,0,\mp 1), \ (0, \pm 1,\mp 1)},
\]
are solutions of \eqref{eq:GFE} irrespective of the signature
$(q,p,r)$. In \cite[Section 4]{PacettiVillagra}, it is proved that the
second pair of solutions induces a singular motive, so it does not
correspond to a newform, whereas the first and last pairs are related
to newforms  with complex
multiplication. This leads to the following question.

\begin{ProblemA}
	\label{question}
Is it true that if $(a,b,c)$ is a primitive
non-trivial solution to (\ref{eq:GFE}) and $f$ is a newform with complex multiplication
then $\overline{\rho}_\id{p}\not\simeq\overline{\rho}_{f,\id{p}}$?
\label{problem}
\end{ProblemA}

When working with Frey hyperelliptic  curves, it is usually not hard to prove that the curve constructed from a non-trivial primitive solution has at least one prime of potentially multiplicative reduction, and, consequently, it cannot have complex multiplication by \cite[Theorem 6]{SerreTate}. Therefore, a positive answer to the posted Question follows from
Conjecture~\ref{conj:large-image}, when $p$ is large enough. 
We prove a similar  result  even when no Frey hyperelliptic curve is available.

\begin{thmA}\label{thm: non-CM}
	Let $q>r>2$ and $p>3$ be prime numbers, with $q$ and $r$ fixed.  Let $(a,b,c)$ be a hypothetical non-trivial primitive solution to
        \begin{equation*}
      	 x^q+y^p+z^r=0,          
        \end{equation*}
	and let $\rho_{\id{p}}$ be the Galois representation attached
        to $(a,b,c)$ in Theorem \ref{thm:rep-existence}. Assume that
        Conjecture \ref{conj:large-image} holds. Then there exists a
        constant $C(q,r)$ (depending only on $q$ and $r$) such that if
        $p>C(q,r)$, for any newform $f$ with complex multiplication we have
        that
        $\overline{\rho}_\id{p}\not\simeq\overline{\rho}_{f,\id{p}}$.
\end{thmA}
%
In~\cite{PacettiVillagra} it was explained how a solution $(a,c,s,\ell,m,n,u,v) \in \Z^2 \times \Z_{\ge 0}^6$ to the equation
\begin{equation}
  \label{eq:ghost}
q^s r^\ell a^q \pm q^m r^n+ q^ur^vc^r=0
\end{equation}
gives rise to a Galois representation with properties that mimic the ones coming from solutions to (\ref{eq:GFE}), making the modular method to fail, even if solutions of \eqref{eq:ghost} do not necessarily correspond to solutions of \eqref{eq:GFE}. 
If all modular forms bypassing the elimination process of the modular method happen to correspond to solutions of~\eqref{eq:ghost} or to have complex multiplication (as usually occurs in the literature), then Theorem~\ref{thm: non-CM}
implies that any solution to~(\ref{eq:GFE}) comes from a solution to the simpler equation~\eqref{eq:ghost}, which is akin to a Thue-Mahler  equation.

A key ingredient for the proof of Theorem \ref{thm: non-CM} is the following result on
  hypergeometric motives (see \S \ref{sec:set-up-motives} for details). Let $q$ and $r$ be prime
  numbers such that $q> r>2$ and let $F=\Q(\zeta_{qr})$.
  \begin{thmB} Let $\z-spec \in \Q \setminus\{0,1\}$. If $v_r(t-1) > 2+\frac{2}{r-1}$, then the
    hypergeometric motive
    $\HGM((\frac{1}{r},-\frac{1}{r}),$ $(\frac{1}{q},-\frac{1}{q})|\z-spec)$ has
    multiplicative reduction at prime ideals of $F$ dividing $r$.
  \end{thmB}
  The result is of particular interest, since it involves studying the
  action of inertia on cyclic covers of $\PP^1$ at primes dividing
  the degree.  Its proof follows the strategy described in \cite{TimDuke}.
  
  \subsection*{Organization of the paper} In
  Section~\ref{Sec:prelims}, we introduce the tools needed to prove
  Theorem~\ref{thm: non-CM}, including results on hypergeometric
  motives and models of curves over discrete valuation rings. In
  Section~\ref{Sec:reduction}, we prove the key result,
  Theorem~\ref{prop:semistable-reduction}, from which
  Theorem~\ref{thm:B} follows. Finally, in Section~\ref{Sec:proof}, we
  prove Theorem~\ref{thm: non-CM}, thereby giving a positive answer to the Question on page \pageref{question}. We want to stress that there is a subtlety while verifying the hypothesis of
  Conjecture~\ref{conj:large-image}, which has not been previously addressed in the literature. In this section,
  we explain the issue and show how to solve it.
	
  \subsection*{Acknowledgments} We thank Tim Dokchitser for explaining
  to us how the results of \cite{TimDuke} can be used even when the
  $\Delta_v$-regularity condition does not hold. We also thank Enric
  Florit for useful discussions regarding building blocks of abelian
  varieties. The first author would also like to thank the second
  author for his great hospitality during his visit to Aveiro.

\section{Preliminaries} 
\label{Sec:prelims}

We present the two main technical tools that will be needed throughout
this paper: hypergeometric motives and models of curves over discrete
valuation rings. We will just recall the details that we need and we
do not intend to provide a full reference for either object. Instead,
we refer the interested reader to \cite[Section 2]{GP} and
\cite{TimDuke} for references.

\subsection{Hypergeometric motives}
\label{sec:set-up-motives}

In this section, we recall the main properties of the rank-$2$
hypergeometric motive
$\HGM((\frac{1}{r},-\frac{1}{r}),(\frac{1}{q},-\frac{1}{q})|t)$ used
to construct the representation $\rho_{\id{p}}$ of
Theorem~\ref{thm:rep-existence}.

Let $q \neq r$ be two rational odd primes. Let $F=\Q(\zeta_{qr})$ and
$K = \Q(\zeta_q)^+\cdot \Q(\zeta_r)^+$ be the compositum of $\Q(\zeta_q)^+$ with $\Q(\zeta_r)^+$. Let $\z-spec \in \Q\setminus\{0,1\}$. Then there exists 
a rank-$2$ hypergeometric motive
$\HGM((\frac{1}{r},-\frac{1}{r}),(\frac{1}{q},-\frac{1}{q})|t)$
defined over $K$. Its set of primes of bad reduction is contained in
$S_{\text{pw}}\cup S_{\text{pt}}$, where
  \[
    S_{\text{pw}}:=\{\lambda \subset K \; : \; \lambda  \mid qr\},
  \]
  is the set of potentially wild primes and
  \[
    S_{\text{pt}}:=\{\lambda  \subset K \; : \; \lambda  \nmid qr \text{ and }v_{\lambda }(\z-spec (\z-spec-1))\neq 0\},
  \]
  is the set of potentially tame primes.
At the potentially tame primes, the image of inertia can be easily computed (see \cite[Corollary 6.14]{GPR}).
Our goal is to describe the behaviour at the potentially wild primes when $v_q(t-1)$ (respectively $v_r(t-1)$) is sufficiently large.

We note that the aforementioned motive can be realized by using the middle cohomology of an Euler curve. Concretely, set 
\begin{equation}
  \label{eq:parameters}
  A:=q-r, \quad B:=-q-r, \quad C:=q+r, \quad D:=-r, \quad N:=qr.
\end{equation}
The \textit{Euler curve} is defined by
\begin{equation} \label{eq:euler-curve-generic}
\bfC : y^{N} = x^{A}(1-x)^{B} (1-\z-spec x)^C\z-spec^{D}.    
\end{equation}
Over $F$, the multiplicative group $\mubb$ of $N$-th roots of unity
acts on $\bfC$ via $\zeta_{N}\cdot(x,y)=(x,\zeta_{N}y)$. The de Rham cohomology
group~$H^1(\bfC,F):=H^1(\bfC,\Q)\otimes_{\Q}F$ admits a decomposition 
\begin{equation}
  \label{eq:decomposition}
  H^1(\bfC,F) = \bigoplus_{\chi} H^1_\chi(\bfC,F),
\end{equation}
where the sum ranges over characters of $\mubb$ and $H^1_\chi(\bfC,F)$
denotes the $\chi$-isotypical component. 
Let
$\chi_0:\mubb\rightarrow F^\times$ denote the identity character.  Up
to a Tate twist, the hypergeometric motive $\HGM((\frac{1}{r},-\frac{1}{r}),(\frac{1}{q},-\frac{1}{q})|\z-spec)$ is
defined as the $\chi_0$-isotypical component of \eqref{eq:decomposition} (see \cite[Definition
6.10]{GPR} for extra details).
To define the motive over $K$ rather than $F$, one needs to consider quotients of
$\bfC$ given by some involutions. However, for our purposes (to
prove that the reduction type of the motive is potentially
multiplicative) it is enough to work with the curve $\bfC$ itself.

\subsection{Models of curves over DVRs} \label{sec:models-over-DVRs}
As seen in the previous section, we need to study the curve
$\bfC$, which is a singular non-hyperelliptic curve. To do so,
we shall apply the techniques developed in \cite{TimDuke}.  We begin with a brief overview of the results and definitions of Dokchitser's
article. We refer the reader to \S3 of loc.\ cit.\ for further
details.

Let $K$ be a discretely valued field with valuation
$ v$, ring of integers $\Om_K$, uniformizer
$\pi$, and residue field $k$.  Consider a bivariate polynomial
\[
  P=\sum_{i,j}a_{i,j}x^iy^j\in K[x,y].
\]
Its Newton polygon is defined by the convex hull of the points $(i,j)$, that is
\[
  \Delta=\operatorname{conv}\{(i,j)\in\Z^2\mid
  a_{i,j}\neq0\}\subset\RR^2.
\]
Similarly, define the Newton polytope of $P$ with respect to the valuation $v$ by
\[
  \Delta_v= \operatorname{lch} \{(i,j,v(a_{i,j}))\mid a_{i,j}\neq0\}
  \subset\RR^2 \times \RR,
\]
where $\operatorname{lch}$ stands for the lower convex hull, i.e., the points $(a_0, b_0, c_0)$ in the convex hull such that $(a_0, b_0,c_0)-(0,0,\varepsilon)$ does not lie in the convex hull for any $\varepsilon > 0$. Finally, we note that 
the projection $p:\Delta_v\to \Delta$ is a homeomorphism
onto its image, and thus the valuation induces a piecewise affine
function $v:\Delta\rightarrow\RR$.
	
\begin{defi}
	\label{defi:vfacesbuena}
  The images in $\Delta$ of the open faces of $\Delta_v$ under the
  above projection are called the \emph{$v$-vertices},
  \emph{$v$-edges}, and \emph{$v$-faces}, according to whether the
  corresponding face of $\Delta_v$ has dimension $0$, $1$, or $2$.
  Thus every $v$-vertex is a lattice point, every $v$-edge (denoted by
  $L$) is homeomorphic to an open interval, and every $v$-face
  (denoted by $F$) is homeomorphic to an open disc. 
\end{defi}
For any $v$-edge $L$ and $v$-face $F$, define
\[
  L(\Z)=L\cap\Z^2, \qquad F(\Z)=F\cap\Z^2.
\]
Similarly, $\Delta(\Z)=\operatorname{int}(\Delta)\cap\Z^2$.
The notation $\bar L(\Z)$, $\bar F(\Z)$, and $\bar\Delta(\Z)$ refers to the previous objects together with the
boundary points.
The subscript $\Z$ on a set indicates the subset of elements satisfying the additional condition that $v(P) \in \Z$. For instance,
\[
   F(\Z)_{\Z} = \{P\in F(\Z)\mid v(P)\in\Z\}.
\]
Finally, for a
$v$-edge or $v$-face $\lambda$, we denote by $\delta_\lambda$ the
common denominator of the numbers $v(P)$ with $P\in\bar\lambda(\Z)$.
	
\begin{defi}[Restriction] \label{def:restriction} Let
  $g=\sum_{\mathbf{i}\in\Z^n} c_{\mathbf{i}}\mathbf{x}^{\mathbf{i}}
  \in K[\mathbf{x}]$ where $\mathbf{x}=(x_1,\dots,x_n)$, and let
  $\emptyset\neq S\subset\Z^n$.  Let $\Lambda$ be the smallest affine
  lattice satisfying $ S\subset\Lambda\subset\Z^n$, and suppose that
  $\Lambda$ has rank $r$.  Choose an affine lattice isomorphism
  $ \phi:\Z^r\longrightarrow\Lambda$, write
  $\mathbf{y}=(y_1,\dots,y_r)$, and define
		\[
		g|_S
		=
		\sum_{\mathbf{i}\in\phi^{-1}(S)}
		c_{\phi(\mathbf{i})}\mathbf{y}^{\mathbf{i}}.
		\]
		Different choices of $\phi$ modify $g|_S$
		only by a monomial factor and a
		$\GL_r(\Z)$-change of variables.		
		For a $v$-vertex, $v$-edge, or $v$-face $\lambda$, we set
		\[
		P|_\lambda
		=
		P|_{\bar\lambda(\Z)_{\Z}}.
		\]
	\end{defi}
	
	\begin{defi}[Reduction]\label{def:reduction}
		Let $h\in K[x,y]$.
		Assume that there exist integers $c,m,n$ such that
		$\tilde h(x,y)
		=
		\pi^c h(\pi^m x,\pi^n y)$
		has coefficients in $\Om_K$, and that
		$\tilde h\bmod\pi$
		has the same Newton polygon as $h$.
		The polynomial $
		\tilde h\bmod\pi\in k[x,y]
		$
		is then called the \emph{reduction} of $h$, and will be denoted by
		$\bar h$.
\end{defi}
	
Given a $v$-edge or $v$-face $\lambda$, define the $k$-scheme
\begin{equation}
  \label{eq:reduction}
X_\lambda \colon\{\overline{P|_\lambda}=0\}.   
\end{equation}
\begin{defi}\label{defi:v-faces}
A $v$-face $F$ is \emph{$\Delta_v$-regular} if the variety $X_F$ is smooth over the residue field $k$. The same definition applies to edges $L$ that are not on the boundary of $\Delta$.
\end{defi}
The definition at edges lying on the boundary is a little more
involved (see \cite[Definition 3.9]{TimDuke}). The zero locus of $P$ is said to be $\Delta_v$-regular if all its faces and edges are $\Delta_v$-regular.

The model of the Euler curve that we shall consider in Section \ref{Sec:reduction} is not $\Delta_v$-regular, but, as explained in \S 10 of \cite{TimDuke}, having $\Delta_v$-regular faces and edges is sometimes enough to apply \cite[Theorem 3.14]{TimDuke}. By gluing the desingularizations, this allows to determine the semistable reduction of $\bfC$ (see Remark 10.1 of loc. cit.).
\section{\texorpdfstring{Reduction at potentially wild primes}{Reduction at wild primes of our motive}}
  \label{Sec:reduction}
Given    two  distinct  odd  prime  numbers $q,r$ and $t\in\Q\setminus\{0,1\}$ we consider  the hypergeometric motive
\begin{equation} \label{eq:motive-qr}
	\HGM\left(\left(\frac{1}{r},-\frac{1}{r}\right),\left(\frac{1}{q},-\frac{1}{q}\right)\Big|\z-spec\right).
\end{equation}
Let $(a,b,c)$ be a hypothetical solution to (\ref{eq:GFE}). The specialization of (\ref{eq:motive-qr}) at $t=-\frac{a^q}{c^r}$ 
gives rise to the representation $\rho_{\id{p}}$ in Theorem \ref{thm:rep-existence} \cite{GP}.

  It  follows  from
\cite[Proposition 4.2]{GPR} that for any  pair of integers $m,n$ prime
to $qr$, the motive in (\ref{eq:motive-qr}) and the motive
\begin{equation} \label{eq:motive2} \HGM\left(\left(\frac{n}{r},-\frac{n}{r}\right),\left(\frac{m}{q},-\frac{m}{q}\right)\Big|\z-spec\right)
\end{equation}
are Galois conjugate. Then it is enough to study the case where $m$
and $n$ satisfy that $nq-mr=1$, so the Euler curve associated to the motive in (\ref{eq:motive2}) has equation
\begin{equation}
  \label{eqn:Euler}
  \bfC : y^{qr} = x(1-x)^B (1-\z-spec x)^C\z-spec^D,
\end{equation}
where $B=-mr-nq$, $C=mr+nq$ and $D=-mr$.
After a simple change of variables, we may furthermore assume that
$0 < B, C < qr$ and that $B+C = qr$.

Let $\zeta_n$ be a primitive $n$-th root of unity and let
$F=\Q(\zeta_{qr})$. Fix a prime ideal $\id{r}$ of $F$ lying above $r$,
and write $F_{\id{r}}$ for the completion of $F$ at $\id{r}$. Let
$\pi$ be a local uniformizer. We denote by $v_{\id{r}}$ the normalized
valuation at $\id{r}$ over $F$ and by $v_r$ the valuation at $r$ over $\Q$.

\begin{thm}
  \label{prop:semistable-reduction}
  Let $\z-spec \in \Q \setminus\{0,1\}$ satisfy $v_{r}(\z-spec-1) > 2+\frac{2}{r-1}$. Then the curve
  $\bfC$ has a semistable model over $F_\id{r}$ whose toric part has full
  rank.
\end{thm}

\begin{proof}
Set $k := v_{\id{r}}(t-1)=(r-1)v_r(t-1) > 2r$.   Since $\z-spec$ is a unit, we can make a change of variables and
  remove the factor $\z-spec^D$ in (\ref{eqn:Euler}) (because the
  extension $F_\id{r}(\sqrt[q]{\z-spec})/F_\id{r}$ is unramified at $\id{r}$).
  Translate the singular point $(1,0)$ to the origin (using the change
  of variables $x\mapsto x+1, \ y\mapsto -y$), so the new
  equation reads $P(x,y)=0$ for 
\begin{equation}
  \label{eqn:uglyfxy}
  P(x,y)=
  y^{qr}
  +(-x)^B\sum_{i=0}^{C+1}(1-t)^{C-i}(-t)^{i-1}\left(\binom{C}{i-1}-t\binom{C+1}{i}\right)x^i.
\end{equation}
Keeping the terminology of Section \ref{sec:models-over-DVRs}, let $\Delta$ denote the Newton polygon of $P(x,y)$.
Let $L_1$ be the line connecting the points $(qr,0)$ and $(0,qr)$, with
equation
\begin{equation}
  \label{eq:1}
  L_1:x+y = qr.
\end{equation}
This line splits the Newton polygon $\Delta$   as $\Delta=F_1 \cup F_2$, where $F_1$ is determined by the points where $x+y \le qr$ and  $F_2$ by the
points where $x+y \ge qr$. Let $v_1:F_1\to\RR$ be the affine function defined  by
\begin{equation}
	\label{eq:affine1}
	v_1(x,y) = qrk-k(x+y),
\end{equation} 
and $v_2:F_2\to\RR$ be the null function. 

\vspace{5pt}

\noindent{\bf Claim 1:} $F_1$ and $F_2$ are the $v$-faces of $\Delta$,
as defined in Definition \ref{defi:vfacesbuena} (see
Figure~\ref{fig:Newton-polytope}\footnote{ Large numbers indicate
  that the corresponding coefficient of $P$ has exactly the specified valuation.}).

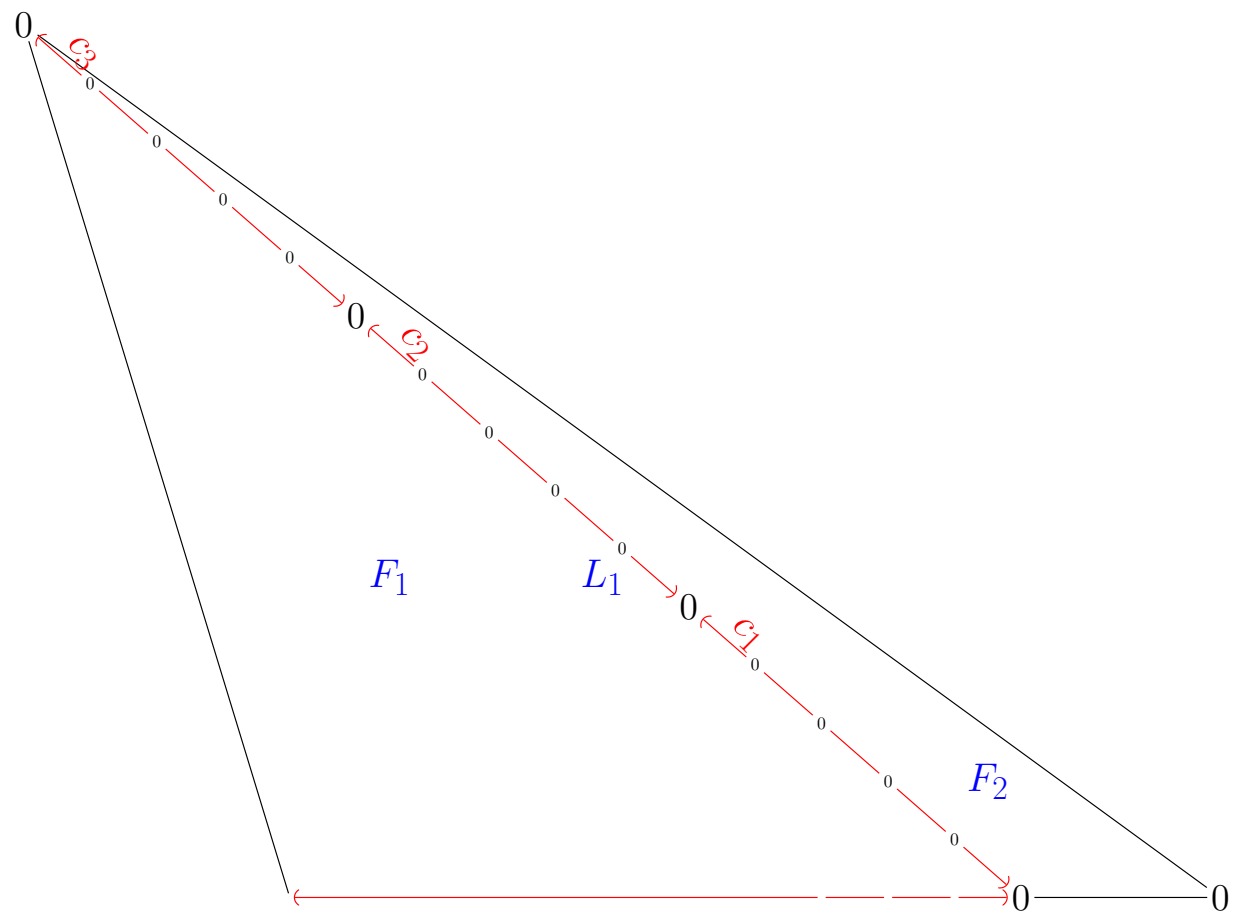
\begin{figure} 
	\centering
	\begingroup
	\resizebox{0.35\textwidth}{!}{
		\def\eqn{}
		\def\p{5}
		\def\deltav{\noindent\pbox{10cm}{\begin{tikzpicture}[xscale=0.8,yscale=0.7,
					fname/.style={blue,scale=0.55},
					lname/.style={scale=0.55,sloped,red,above=-0.07em,near end},
					lrg/.style={scale=0.9,inner sep=0.2em},
					sml/.style={scale=0.55},
					lin/.style={-,shorten <=-0.07em,shorten >=-0.07em}]
					\node[fname] at (5.5,5.5) {\Huge$F_1$};
					\node[fname] at (14.5,2) {\Huge$F_2$};
					\node[fname] at (8.7,5.5) {\Huge$L_1$};
					\node[lrg] at (0,15) (1) {\Large 0};
					\node[sml] at (1,14) (2) {0};
					\node[sml] at (2,13) (3) {0};
					\node[sml] at (3,11) (4) {};
					\node[sml] at (3,12) (5) {0};
					\node[sml] at (4,10) (6) {};
					\node[sml] at (4,11) (7) {0};
					\node[sml] at (5,9) (8) {};
					\node[lrg] at (5,10) (9) {\Large 0};
					\node[sml] at (6,7) (10) {};
					\node[sml] at (6,8) (11) {};
					\node[sml] at (6,9) (12) {0};
					\node[sml] at (7,6) (13) {};
					\node[sml] at (7,7) (14) {};
					\node[sml] at (7,8) (15) {0};
					\node[sml] at (8,5) (16) {};
					\node[sml] at (8,6) (17) {};
					\node[sml] at (8,7) (18) {0};
					\node[sml] at (9,3) (19) {};
					\node[sml] at (9,4) (20) {};
					\node[sml] at (9,5) (21) {};
					\node[sml] at (9,6) (22) {0};
					\node[sml] at (10,2) (23) {};
					\node[sml] at (10,3) (24) {};
					\node[sml] at (10,4) (25) {};
					\node[lrg] at (10,5) (26) {\Large 0};
					\node[lrg] at (4,0) (27) {};
					\node[sml] at (11,1) (28) {};
					\node[sml] at (11,2) (29) {};
					\node[sml] at (11,3) (30) {};
					\node[sml] at (11,4) (31) {0};
					\node[lrg] at (12,0) (32) {};
					\node[sml] at (12,1) (33) {};
					\node[sml] at (12,2) (34) {};
					\node[sml] at (12,3) (35) {0};
					\node[lrg] at (13,0) (36) {};
					\node[sml] at (13,1) (37) {};
					\node[sml] at (13,2) (38) {0};
					\node[lrg] at (14,0) (39) {};
					\node[sml] at (14,1) (40) {0};
					\node[lrg] at (15,0) (41) {\Large 0};
					\node[lrg] at (18,0) (42) {\Large 0};
					\draw[lin]
					(1) edge (42)
					(1) edge[red,(-] node[lname] {\Huge$c_3$} (2) (2) edge[red] (3) (3) edge[red] (5) (5) 
					edge[red] (7) (7) edge[red,-)] (9) (9) edge[red,(-] node[lname] {\Huge$c_2$} (12) (12) 
					edge[red] (15) (15) edge[red] (18) (18) edge[red] (22) (22) edge[red,-)] (26) 
					(26) edge[red,(-] node[lname] {\Huge$c_1$} (31) (31) edge[red] (35) (35) edge[red] (38) 
					(38) edge[red] (40) (40) edge[red,-)] (41)
					(41) edge (42)
					(27) edge[red,(-] node[lname] {} (32) (32) edge[red] (36) (36) edge[red] (39) 
					(39) edge[red,-)] (41)
					(1) edge (27)
					;
			\end{tikzpicture}}\hfill}
		\deltav}
	\endgroup
	\caption{Newton polytope of $P(x,y)$ for $r=5$ and $q=3$.}
	\label{fig:Newton-polytope}
\end{figure}

\vspace{2pt}

\noindent\textit{Proof of Claim 1.} Write $
  P(x,y)=\sum_{i,j} a_{i,j}x^iy^j.$
The  non-zero coefficients of $P(x,y)$ are  $a_{0,qr}$ and $a_{i,0}$, for $i=B, \dots, qr+1$.
 Then, Claim 1 follows easily from the data in Table~\ref{table:val-1}, which contains information on the
 valuation of the non-zero coefficients $a_{i,j}$,  together with the fact that $B+C=qr$.
\begin{table} 
\begin{tabular}{|c|c|c||c|c|c|}
  \hline
 $i$ & $j$ & $v_\id{r}(a_{i,j})$  & $i$ & $j$ & $v_\id{r}(a_{i,j})$ \\
  \hline \hline
  $0$ & $qr$ & $0$ & $qr-1$ & $0$ & $k$ \\
  $qr$ & $0$ & $0$ & $qr+1$ & $0$ & $0$ \\
  $B$ & $0$ & $Ck$ & $[B, qr-2]$ & 0 & $\ge(qr-i)k$ \\
  \hline 
\end{tabular}
\caption{Valuations of the coefficients appearing in (\ref{eqn:uglyfxy}).}
\label{table:val-1}
\end{table}

\vspace{5pt}

\noindent {\bf Claim 2:} The two faces $F_1$ and $F_2$ are
$\Delta_v$-regular.

\vspace{2pt}

\noindent \textit{Proof of Claim 2.} By Definition \ref{def:restriction}, we have
\[
{P|_{F_1}}=y^{qr}
+(-x)^B\sum_{i=0}^{C}(1-t)^{C-i}(-t)^{i-1}\left(\binom{C}{i-1}-t\binom{C+1}{i}\right)x^i
\]
and
\[
P|_{F_2}=t^{C-1}(C-t(C+1))x^{qr}-t^Cx^{qr+1}+y^{qr}.
\]
Let $u \equiv \frac{1-t}{\pi^k} \pmod{\id{r}}$. Then the reduction $\overline{P|_{F_1}}$ equals
\[
\overline{P|_{F_1}}=y^{qr} + u^C(-1)^Bx^B \sum_{i=0}^C \binom{C}{i}(-x/u)^i = y^{qr}+(-x)^B(u-x)^C,
\]
which defines a smooth curve in the  $\Gm_m^2$, where $\Gm_m$ stands for the multiplicative group.
A similar computation, using the fact that $C\equiv 1 \pmod r$,  proves the result for $F_2$.

\vspace{5pt}

If $P_1,P_2$ are points in $\Z^2$, we use the notation $L:=P_1 \vdashv P_2$ to
define the segment that joins $P_1$ and $P_2$. Consider  the following edges
(appearing in Figure~\ref{fig:Newton-polytope}).
\begin{align*}
  L_2&:= (0,qr) \vdashv(qr+1,0), &L_3&:=(0,qr) \vdashv (B,0),\\
  L_4&:=(B,0) \vdashv(qr,0), & L_5&:=(qr,0)\vdashv (qr+1,0),
\end{align*}
along with the edge $L_1 = (qr,0)\vdashv(0, qr)$ defined in \eqref{eq:1}. Then the  reductions of $P$ at the edges are
\begin{equation}
  \label{eqn:schemeL1}
  \overline{P|_{L_1}}
  =T^{qr}-1
  =\prod_{j=1}^{q}\bigl(T-\zeta_q^j\bigr)^r, \quad \text{and} \quad
\overline{P|_{L_4}}=(-T)^B(u-T)^C,
 \end{equation}
corresponding to edges which are not $\Delta_v$-regular and 
\[
\overline{P|_{L_2}}=1- T, \quad \overline{P|_{L_3}}=1-T, \quad \overline{P|_{L_5}}=-1-T,
\]
which are $\Delta_v$-regular. 
The edge $L_1$
 lies on the boundary of $F_1$ and $F_2$ and, by \eqref{eqn:schemeL1}, it consists of $q$ singular points. Then, 
the graph of components is illustrated in
Figure~\ref{fig:componentes}, where the red edges mean that the resulting model is not $\Delta_v$-regular. 
\begin{figure} 
\centering
\begingroup
\resizebox{0.5\textwidth}{!}{%
\def\graph{\pbox{10cm}{\begin{tikzpicture}[xscale=1,yscale=0.9]
\draw[thick,shorten >=2] (-1/3,0)--(86/15,0) node[inner sep=2pt,blue,above left,scale=0.8] {1} node[inner sep=2pt,below left,scale=0.8] {$F_1$};
\node[red,label={[above,inner sep=0pt,scale=0.7,red]}] at (0,0) {$\bullet$};
\draw[shorten <=-3pt, shorten >=-3pt, thick, red] (7/10,0)--node[inner 	sep=2pt,scale=0.8,blue,right] {$c_1$} (7/10,4/5);
\draw[shorten <=-3pt, shorten >=-3pt, thick, red] (3/2,0)--node[inner 	sep=2pt,scale=0.8,blue,right] {$c_2$} (3/2,4/5);
\draw --node[inner sep=2pt,scale=0.8,blue,right] {$ \cdots$} (50/10,4/5);
\draw[shorten <=-3pt, shorten >=-3pt, thick, red] (39/10,0)--node[inner sep=2pt,scale=0.8,blue,right] {$c_q$} (39/10,4/5);
\draw[thick,shorten >=2] (11/30,4/5)--(86/15,4/5) node[inner sep=2pt,blue,above left,scale=0.8] {1} node[inner sep=2pt,below left,scale=0.8] {$F_2$};
\end{tikzpicture}}\hfill}
\graph}
\endgroup
\caption{Graph of components of the Euler curve $\bfC: P(x,y)=0$.}
\label{fig:componentes}
\end{figure}
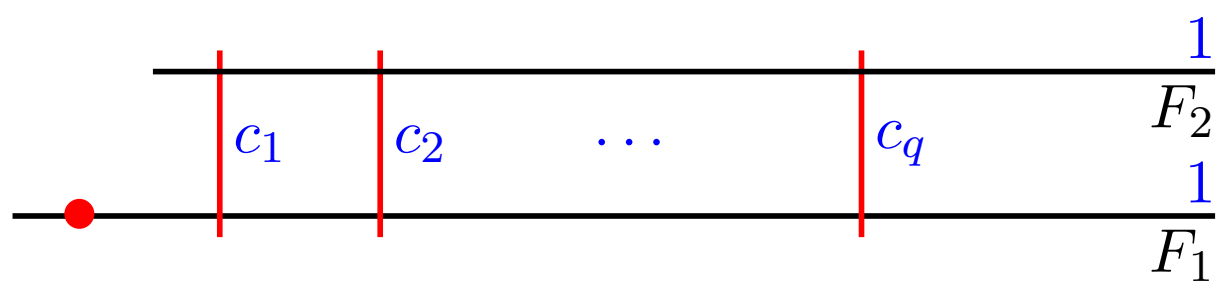

\vspace{5pt}
	
Fix $j \in \{1,\ldots, q\}$. To resolve the line $c_j$ in
Figure~\ref{fig:componentes} we need to  find a model that satisfies
the hypothesis of \cite[Theorem 3.14]{TimDuke}. To do so, we translate the
singular point $\zeta_q^j$ of (\ref{eqn:schemeL1}) to a point not
lying in $\Gm_m$.
The change of variable $y\to y+\zeta_q^jx$
transforms equation~(\ref{eqn:uglyfxy}) into the model $Q(x,y)=0$, where
\[
  Q(x,y)=  (y+\zeta_q^jx)^{qr}
  +(-x)^B\sum_{i=0}^{C+1}(1-t)^{C-i}(-t)^{i-1}\left(\binom{C}{i-1}-t\binom{C+1}{i}\right)x^i.
\]
Define the following segments:
\begin{align*}
  L_1'&:= (0,qr) \vdashv (qr-r,r), &  L_2'&:=(qr-r,r) \vdashv (qr-1,1),\\
  L_3'&:=(qr-1,1) \vdashv (qr-1,0),&  L_4'&:=(qr-1,1) \vdashv(qr+1,0),\\
  L_5'&:= (0,qr) \vdashv (B,0),&  L_6'&:=(0,qr) \vdashv (qr+1,0),\\
  L_7'&:=(B,0) \vdashv(qr-1,0),&  L_8'&:=(qr-r,r) \vdashv(qr-1,0),\\
 L_9'&:=(qr-r,r) \vdashv (qr+1,0),&  L_{10}'&:=(qr-1,0) \vdashv(qr+1,0).
\end{align*}
The segments decompose the polygon $\Delta$ in the form
\[
  \Delta = F_1'\cup F_2'\cup F_3'\cup F_4'\cup F_5',
\]
where
\begin{itemize}
\item $F_1'$ is the polygon with boundary $L_3' \cup L_4'\cup L_{10}'$.
\item $F_2'$ is the polygon with boundary $L_2' \cup L_4' \cup L_9'$.
\item $F_3'$ is the polygon with boundary $L_1' \cup L_8' \cup L_7' \cup L_5'$.
\item $F_4'$ is the polygon with boundary $L_2' \cup L_3' \cup L_8'$.
\item $F_5'$ is the polygon with boundary $L_1' \cup L_9' \cup L_6'$.
\end{itemize}
Define the following affine functions $v'_i:F_i'\to \RR$:
\[
v'_1(x,y)=-\frac{k}{2}x+(r-1-k)y+\frac{k(qr+1)}{2}, \quad v'_2(x,y)=-rx-(r+1)y+r(qr+1),
\]
\[
  v'_3(x,y)=qrk-k(x+y), \quad v'_4(x,y)=(r-k)x+(r-1-k)y+r(qk-qr+1), \quad v'_5(x,y)=0.
\]

\vspace{5pt}

\noindent {\bf Claim 3:} $F_i$'s are the $v$-faces of the Newton polytope of the polynomial $Q(x,y)$ (see Figure~\ref{fig:desing1}).

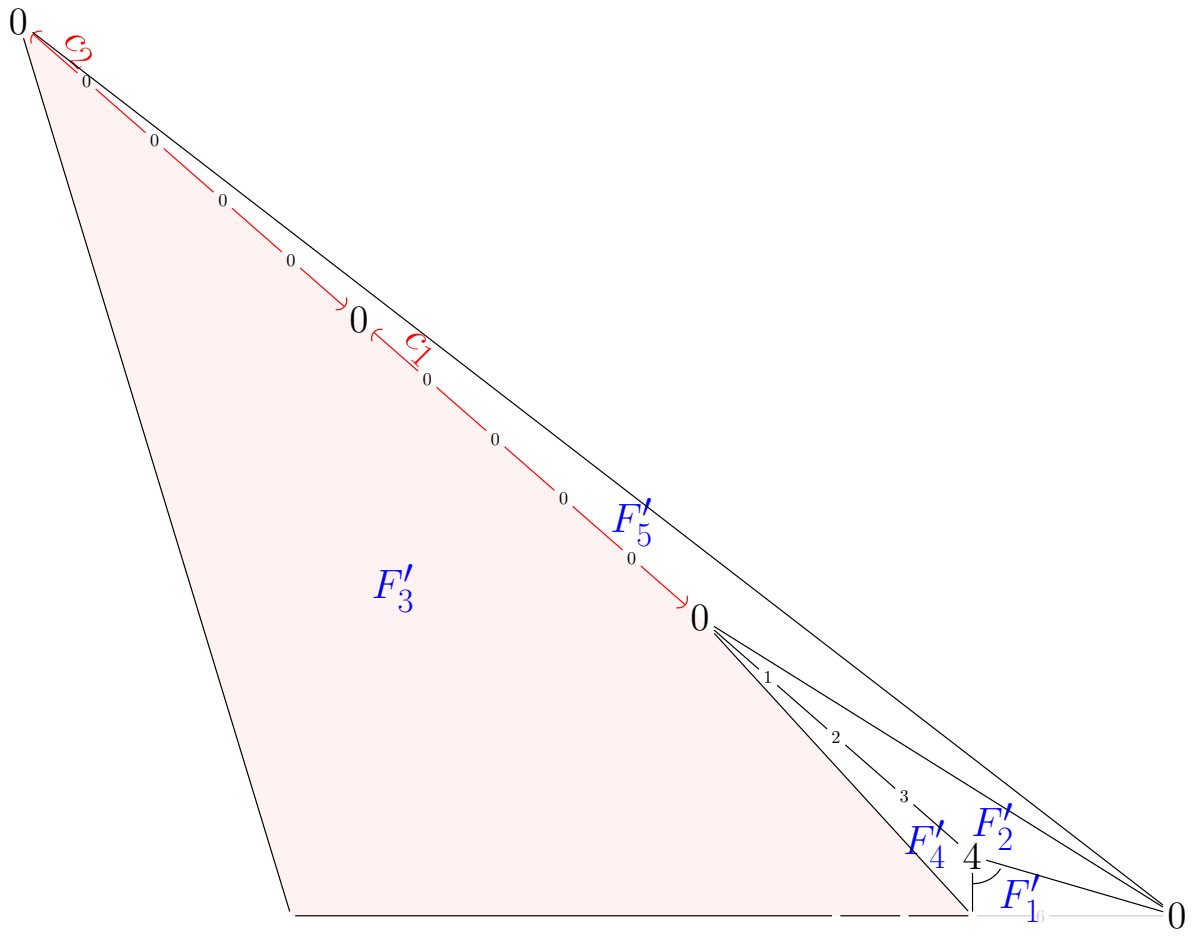
\begin{figure} 
	\centering
	\begingroup
	\resizebox{0.35\textwidth}{!}{%
		\def\p{5}
		\def\deltav{\noindent\pbox{10cm}{\begin{tikzpicture}[xscale=0.8,yscale=0.7,
					fname/.style={blue,scale=0.55},
					lname/.style={scale=0.55,sloped,red,above=-0.07em,near end},
					lrg/.style={scale=0.9,inner sep=0.2em},
					rem/.style={black!20,thin},
					sml/.style={scale=0.55},
					lin/.style={-,shorten <=-0.07em,shorten >=-0.07em}]
					\fill[red!5] (4,0)--(14,0)--(10,5)--(0,15)--cycle;
					\node[fname] at (14.7,0.3) {\Huge $F_1'$};
					\node[fname] at (13.3,1.2) {\Huge $F_4'$};
					\node[fname] at (14.3,1.5) {\Huge $F_2'$};
					\node[fname] at (5.5,5.5) {\Huge $F_3'$};
					\node[fname] at (9,6.6) {\Huge $F_5'$};
					\node[lrg] at (0,15) (1) {\Large 0};
					\node[sml] at (1,14) (2) {0};
					\node[sml] at (2,13) (3) {0};
					\node[sml,red!70!black] at (3,11) (4) {};
					\node[sml] at (3,12) (5) {0};
					\node[sml,red!70!black] at (4,10) (6) {};
					\node[sml] at (4,11) (7) {0};
					\node[sml,red!70!black] at (5,9) (8) {};
					\node[lrg] at (5,10) (9) {\Large 0};
					\node[sml,red!70!black] at (6,7) (10) {};
					\node[sml,red!70!black] at (6,8) (11) {};
					\node[sml] at (6,9) (12) {0};
					\node[sml,red!70!black] at (7,6) (13) {};
					\node[sml,red!70!black] at (7,7) (14) {};
					\node[sml] at (7,8) (15) {0};
					\node[sml,red!70!black] at (8,5) (16) {};
					\node[sml,red!70!black] at (8,6) (17) {};
					\node[sml] at (8,7) (18) {0};
					\node[sml,red!70!black] at (9,3) (19) {};
					\node[sml,red!70!black] at (9,4) (20) {};
					\node[sml,red!70!black] at (9,5) (21) {};
					\node[sml] at (9,6) (22) {0};
					\node[sml,red!70!black] at (10,2) (23) {};
					\node[sml,red!70!black] at (10,3) (24) {};
					\node[sml,red!70!black] at (10,4) (25) {};
					\node[lrg] at (10,5) (26) {\Large 0};
					\node[lrg] at (4,0) (27) {};
					\node[sml,red!70!black] at (11,1) (28) {};
					\node[sml,red!70!black] at (11,2) (29) {};
					\node[sml,red!70!black] at (11,3) (30) {};
					\node[sml] at (11,4) (31) {1};
					\node[lrg] at (12,0) (32) {};
					\node[sml,red!70!black] at (12,1) (33) {};
					\node[sml,red!70!black] at (12,2) (34) {};
					\node[sml] at (12,3) (35) {2};
					\node[lrg] at (13,0) (36) {};
					\node[sml,red!70!black] at (13,1) (37) {};
					\node[sml] at (13,2) (38) {3};
					\node[lrg] at (14,0) (39) {};
					\node[lrg] at (14,1) (40) {\Large4};
					\node[sml,rem] at (15,0) (41) {6};
					\node[lrg] at (17,0) (42) {\Large 0};
					\draw[lin]
					(1) edge (27)
					(26) edge (39)
					(39) edge (40)
					(1) edge[red,(-] node[lname] {\Huge $c_2$} (2) (2) edge[red] (3) (3) edge[red] (5) (5) 
					edge[red] (7) (7) edge[red,-)] (9) (9) edge[red,(-] node[lname] {\Huge $c_1$} (12) (12) 
					edge[red] (15) (15) edge[red] (18) (18) edge[red] (22) (22) edge[red,-)] (26)
					(27) edge (32) (32) edge (36) (36) edge (39)
					(39) edge[rem] (41) (41) edge[rem] (42)
					(40) edge (42)
					(26) edge (31) (31) edge (35) (35) edge (38) (38) edge (40)
					(26) edge (42)
					(1) edge (42)
					;
					\draw [-] (40) ++(270:1.2em) arc (270:333:1.2em);
			\end{tikzpicture}}\hfill}
		\deltav}
	\endgroup
	\caption{Desingularization along one segment for $r=5$ and $q=3$.}
	\label{fig:desing1}
\end{figure}

\vspace{2pt}

\noindent \textit{Proof of Claim 3.} Write $Q(x,y)= \sum b_{i,j}x^iy^j$. The indices $(i,j)$ for which $b_{i,j}\neq 0$  belong either to the $x$-axis or to the line  $x+y=qr$. In Table~\ref{table:val-2}, we summarize information 
on the valuation of the coefficients $b_{i,j}$. Then,  Claim 3 follows from the data. Note that it is here where the hypothesis $k > 2r$ is used. This is in order to
ensure that the point $(qr-1,1,r-1)$ belongs to the lower convex hull and thus the face $F'_1$ is contractible (see \cite[Definition 5.1]{TimDuke}). 

\begin{table} 
	\begin{tabular}{|c|c|c||c|c|c|}
		\hline
		$i$ & $j$ & $v_\id{r}(b_{i,j})$ & $i$ & $j$ & $v_\id{r}(b_{i,j})$ \\
		\hline   \hline
		$0$ & $qr$ & $0$ & $qr-1$ & $1$ & $r-1$ \\
		$qr-1$ & $0$ & $k$ & $qr-r$ & $r$ & $0$ \\
		$B$ & $0$ & $Ck$ & $qr+1$ & $0$ & $0$ \\
		$[B,qr-2]$ & $0$ & $\ge (qr-i)k$ & $[qr-2,qr-r+1]$ & $qr-i$ & $\ge r-1$ \\                  
		\hline
	\end{tabular}
	\caption{Valuations of the coefficients appearing in $Q(x,y)$.}
	\label{table:val-2}
\end{table}

\vspace{5pt}

\noindent {\bf Claim 4:} The faces $F_2'$ and $F_4'$ are $\Delta_v$-regular and the same holds for the line $L_2'$.

\vspace{2pt}

\noindent \textit{Proof of Claim 4.} A direct computation give that the reduced equations are
  \[
\overline{Q|_{F_4'}}=x^{qr-1}y\zeta_q^{-j}q\frac{r}{\pi^{r-1}}+x^{qr-r}y^{r}\zeta_q^{-rj} q + ux^{qr-1},
\]
\[
\overline{Q|_{F_2'}}=x^{qr-1}y\zeta_q^{-j}q\frac{r}{\pi^{r-1}}+x^{qr-r}y^{r}\zeta_q^{-rj} q - x^{qr+1}
\]
and
\[
\overline{Q|_{L_2'}}=q\zeta_q^{-j}\left(\frac{r}{\pi^{r-1}}1+T^{r-1}\zeta_q^{(1-r)j}\right).
\]
All of them are $\Delta_v$-regular. 

\vspace{5pt}

Then Theorem 3.14 in
\cite{TimDuke} implies that the desingularization along $c_j$ consists
of $r-1$ chains of $\PP^1$'s joining the component $F_1$ with the
component $F_2$  as in Figure~\ref{fig:componentes}. The number
of projective lines depends on, following Dokchitser's notation, the values $s_2^{L_2'}$ and $s_4^{L_2'}$ of Definition 3.12 in loc.\ cit. Indeed, take
$P_1=(qr+1,0)$ and $P_0=(qr,0)$. Since $\delta_{L_2'}=1$, we have that
\[
s_2^{L_2'} =v_2(P_1) - v_2(P_0) = -r \quad \text{and}  \quad s_4^{L_2'}=v_4(P_1) - v_4(P_0) = r-k.
\]
The hypothesis  $v_r(t-1)> 2 + \frac{r-1}{2}$ implies that $-r>r-k$, so the number of projective
lines is $k-2r-1$ by \cite[Theorem 3.14]{TimDuke}. The desingularization looks like
Figure~\ref{fig:desing2}, where the dashed blue line connecting $F_2'$
and $F_4'$ comes from the contractible face $F_1'$. 
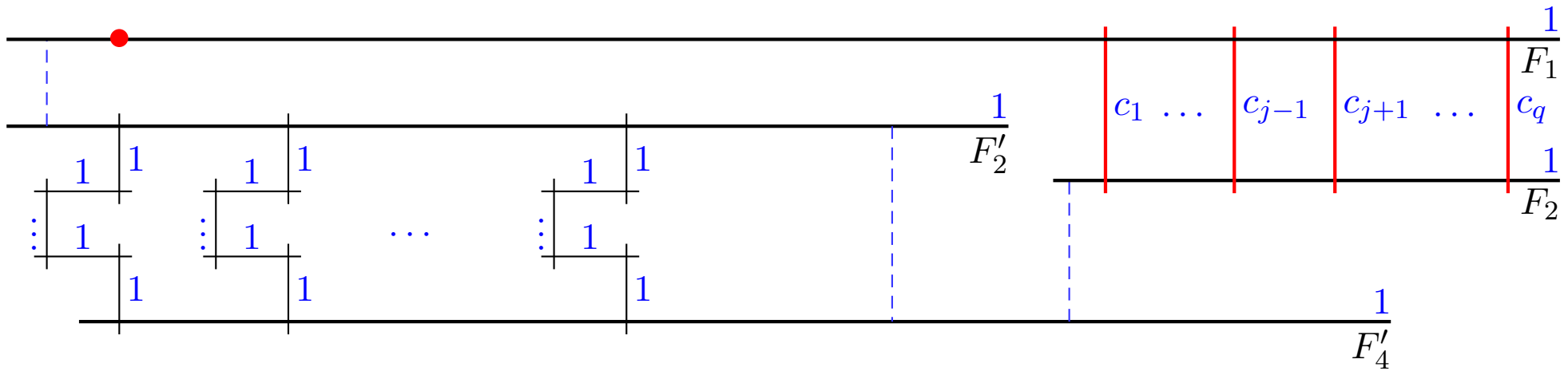
\begin{figure}
	\centering
	\def\graph{\pbox{10cm}{\begin{tikzpicture}[xscale=1,yscale=0.9]
				\draw[thick,shorten >=2] (4/15,0)--(56/5,0) node[inner sep=2pt,blue,above left,scale=0.8] {1} node[inner sep=2pt,below left,scale=0.8] {$F_4'$};
				\draw[shorten <=-3pt, shorten >=-3pt] (3/5,0)--node[inner sep=2pt,scale=0.8,blue,right] {1} (3/5,3/5);
				\draw[shorten <=-3pt, shorten >=-3pt] (3/5,3/5)--node[inner sep=2pt,scale=0.8,blue,above] {1} (0,3/5);
				\draw[shorten <=-3pt, shorten >=-3pt] (0,3/5)--node[inner sep=2pt,scale=0.8,blue,left] {$\vdots$} (0,6/5);
				\draw[shorten <=-3pt, shorten >=-3pt] (0,6/5)--node[inner sep=2pt,scale=0.8,blue,above] {1} (3/5,6/5);
				\draw[shorten <=-3pt, shorten >=-3pt] (3/5,6/5)--node[inner sep=2pt,scale=0.8,blue,right] {1} (3/5,9/5);
				\draw[shorten <=-3pt, shorten >=-3pt] (2,0)--node[inner sep=2pt,scale=0.8,blue,right] {1} (2,3/5);
				\draw[shorten <=-3pt, shorten >=-3pt] (2,3/5)--node[inner sep=2pt,scale=0.8,blue,above] {1} (7/5,3/5);
				\draw[shorten <=-3pt, shorten >=-3pt] (7/5,3/5)--node[inner sep=2pt,scale=0.8,blue,left] {$\vdots$} (7/5,6/5);
				\draw[shorten <=-3pt, shorten >=-3pt] (7/5,6/5)--node[inner sep=2pt,scale=0.8,blue,above] {1} (2,6/5);
				\draw[shorten <=-3pt, shorten >=-3pt] (2,6/5)--node[inner sep=2pt,scale=0.8,blue,right] {1} (2,9/5);
				\draw[shorten <=-3pt, shorten >=-3pt] (24/5,0)--node[inner sep=2pt,scale=0.8,blue,right] {1} (24/5,3/5);
				\draw[shorten <=-3pt, shorten >=-3pt] (24/5,3/5)--node[inner sep=2pt,scale=0.8,blue,above] {1} (21/5,3/5);
				\draw[shorten <=-3pt, shorten >=-3pt] (21/5,3/5)--node[inner sep=2pt,scale=0.8,blue,left] {$\vdots$} (21/5,6/5);
				\draw[shorten <=-3pt, shorten >=-3pt] (21/5,6/5)--node[inner sep=2pt,scale=0.8,blue,above] {1} (24/5,6/5);
				\draw[shorten <=-3pt, shorten >=-3pt] (24/5,6/5)--node[inner sep=2pt,scale=0.8,blue,right] {1} (24/5,9/5);
				\draw--node[inner sep=2pt,scale=0.8,blue,left] {$\cdots$} (33/5,8/5);
				\draw[thick,shorten >=2] (-1/3,9/5)--(241/30,9/5) node[inner sep=2pt,blue,above left,scale=0.8] {1} node[inner sep=2pt,below left,scale=0.8] {$F_2'$};
				\draw[blue!80!white,densely dashed] (7,9/5)--(7,0);
				\draw[blue!80!white,densely dashed] (0,9/5)--(0,13/5);
				\draw[blue!80!white,densely dashed] (254/30,0)--(254/30,13/10);
				\draw[thick,shorten >=2] (250/30,13/10)--(63/5,13/10) node[inner sep=2pt,blue,above left,scale=0.8] {1} node[inner sep=2pt,below left,scale=0.8] {$F_2$};
				\draw[shorten <=-3pt, shorten >=-3pt, thick, red] (263/30,13/10)--node[inner sep=2pt,scale=0.8,blue,right] {$c_1$} (263/30,13/5);
				\draw--node[inner sep=2pt,scale=0.8,blue,right] {$\cdots$} (550/30,19/5);
				\draw[shorten <=-3pt, shorten >=-3pt, thick, red] (295/30,13/10)--node[inner sep=2pt,scale=0.8,blue,right] {$c_{j-1}$} (295/30,13/5);		
					\draw[shorten <=-3pt, shorten >=-3pt, thick, red] (320/30,13/10)--node[inner sep=2pt,scale=0.8,blue,right] {$c_{j+1}$} (320/30,13/5);			
					\draw--node[inner sep=2pt,scale=0.8,blue,right] {$\cdots$} (685/30,19/5);
				\draw[shorten <=-3pt, shorten >=-3pt, thick, red] (363/30,13/10)--node[inner sep=2pt,scale=0.8,blue,right] {$c_{q}$} (363/30,13/5);
				\draw[thick,shorten >=2] (-1/3,13/5)--(63/5,13/5) node[inner sep=2pt,blue,above left,scale=0.8] {1} node[inner sep=2pt,below left,scale=0.8] {$F_1$};
				\node[red,label={[above,inner sep=0pt,scale=0.7,red]}] at (3/5,13/5) {$\bullet$};
	\end{tikzpicture} }}
	\def\charts{}
	\graph\charts
		\caption{Graph of components of the desingularization.}
	\label{fig:desing2}
\end{figure}

The graph of components is obtained by gluing together each of the
desingularizations. The toric rank equals the number of independent
cycles in the dual graph of the graph of components (see \cite[Proposition 3.5]{Grothendieck}). It can
be computed using the well-known formula
\[
  \#E-\#V+1,
\]
where $E$ denotes the set of edges and $V$ the set of vertices (see,
for example, \cite[Theorem 1.9.5(ii)]{graphs}). Note that the formula
is invariant if in Figure~\ref{fig:desing2} we collapse the chain of
$\PP^1$'s into a single one, allowing to consider the simplified dual
graph of Figure~\ref{fig:dualgraph} (the vertices $F_4^{(j)}$ and $F_2^{(j)}$ correspond to the components appearing while taking the desingularization at $c_j$).
\begin{figure} 
	\centering

\begin{tikzpicture}[
	vertex/.style={circle,fill=black,inner sep=2.5pt,outer sep=0pt},
	line width=0.7pt
	]
	
	\draw (0,0) to[bend left=18] (2,2);
	\draw (0,0) to[bend left=18] (2,1);
	\draw (0,0) to[bend right=18] (2,-1);
	\draw (0,0) to[bend right=18] (2,-2);
	
	\draw (2,2) to[bend left=35] (4,2);
	\draw (2,2) to[bend left=20] (4,2);
	\draw (2,2) to[bend right=20] (4,2);
	\draw (2,2) to[bend right=35] (4,2);
	
	\draw (2,1) to[bend left=35] (4,1);
	\draw (2,1) to[bend left=20] (4,1);
	\draw (2,1) to[bend right=20] (4,1);
	\draw (2,1) to[bend right=35] (4,1);
	
	\draw (2,-1) to[bend left=35] (4,-1);
	\draw (2,-1) to[bend left=20] (4,-1);
	\draw (2,-1) to[bend right=20] (4,-1);
	\draw (2,-1) to[bend right=35] (4,-1);
	
	\draw (2,-2) to[bend left=35] (4,-2);
	\draw (2,-2) to[bend left=20] (4,-2);
	\draw (2,-2) to[bend right=20] (4,-2);
	\draw (2,-2) to[bend right=35] (4,-2);
	
	\draw (4,2) to[bend left=18] (6,0);
	\draw (4,1) to[bend left=18] (6,0);
	\draw (4,-1) to[bend right=18] (6,0);
	\draw (4,-2) to[bend right=18] (6,0);
	
	\node[vertex,label=above left:{$F_2$}] at (0,0) {};
	
	\node[vertex,label=above:{$F_4^{(1)}$}] at (2,2) {};
	\node[vertex,label=above:{$F_4^{(2)}$}] at (2,1) {};
	\node[vertex,label=above:{$F_4^{(q-1)}$}] at (2,-1) {};
	\node[vertex,label=above:{$F_4^{(q)}$}] at (2,-2) {};
	
	\node[vertex,label=above:{$F_2^{(1)}$}] at (4,2) {};
	\node[vertex,label=above:{$F_2^{(2)}$}] at (4,1) {};
	\node[vertex,label=above:{$F_2^{(q-1)}$}] at (4,-1) {};
	\node[vertex,label=above:{$F_2^{(q)}$}] at (4,-2) {};
	
	\node[vertex,label=above right:{$F_1$}] at (6,0) {};
	
	\node at (3,2) {$\cdots$};
	\node at (3,2.5) {\footnotesize$r$};
	
	\node at (3,1) {$\cdots$};
	\node at (3,1.5) {\footnotesize$r$};
	
	\node at (3,-1) {$\cdots$};
	\node at (3,-0.5) {\footnotesize$r$};
	
	\node at (3,-2) {$\cdots$};
	\node at (3,-1.5) {\footnotesize$r$};
	
	\node at (3,0.3) {$\vdots$};
	
\end{tikzpicture}
\caption{Simplified dual graph.}
\label{fig:dualgraph}
\end{figure}
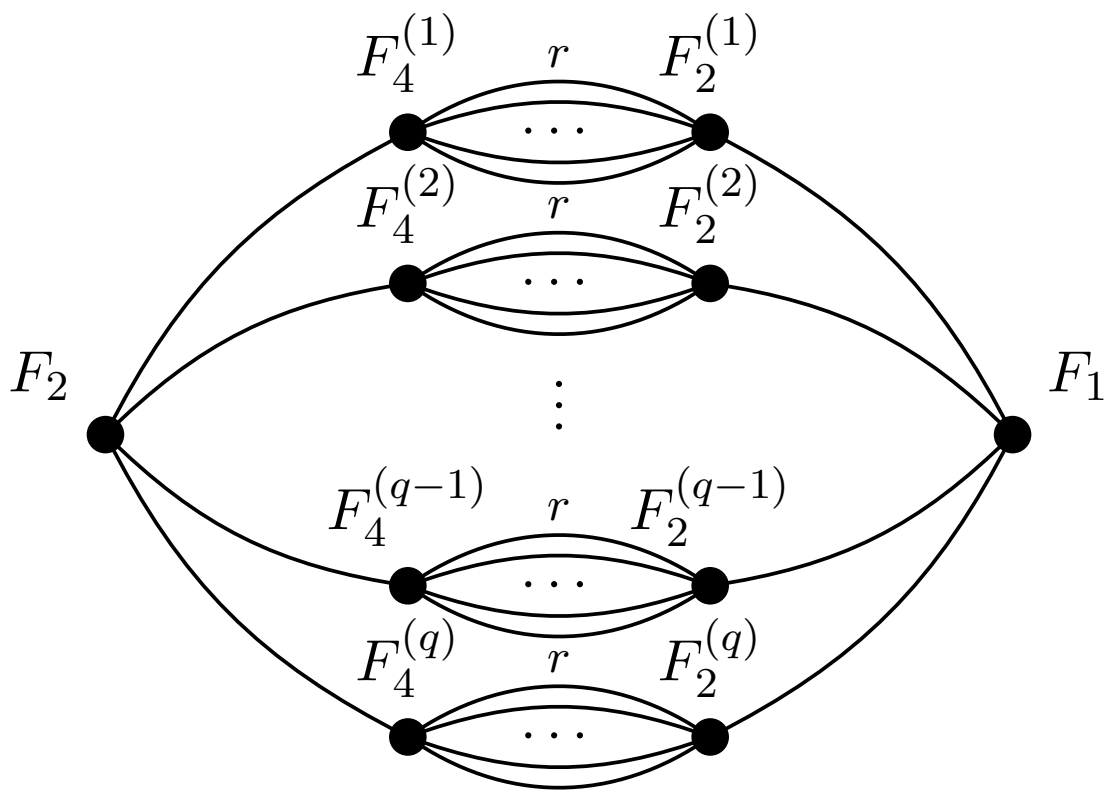

It
readily follows from Figure \ref{fig:dualgraph} that $\#E = qr+2q$,
$\#V = 2q+2$ and so the total of linearly independent cycles is
\[qr+2q-(2q+2)+1 = qr-1,\]
which matches the genus of the curve (see \cite[Theorem~4.1]{MR2005278}).
\end{proof}

  \begin{thmB} \label{thm:B} Let $q, r>2$ be primes. Let $\z-spec \in \Q \setminus\{0,1\}$.  If $v_r(t-1) > 2+\frac{2}{r-1}$, then the
	hypergeometric motive
	$\HGM((\frac{1}{r},-\frac{1}{r}),(\frac{1}{q},-\frac{1}{q})|\z-spec)$ has
	multiplicative reduction at prime ideals of $F$ dividing $r$.
\end{thmB}

\begin{proof}
  From Definition 6.8 in \cite{GPR},  the motive
  $\HGM((\frac{1}{r},-\frac{1}{r}),(\frac{1}{q},-\frac{1}{q})|\z-spec)$
  is a Tate twist of the motive associated to the Euler curve \eqref{eqn:Euler}. The result  then follows from Theorem~\ref{prop:semistable-reduction}.
\end{proof}

\begin{remark}
	\label{rem:simmetry}
	Theorem \ref{thm:B} also holds for the prime $q$. This is due to the fact that, by 
	\cite[Theorem 6.2]{GPR}, we have that
	\[
	\HGM\left(\left(\frac{1}{r},-\frac{1}{r}\right),\left(\frac{1}{q},-\frac{1}{q}\right)\Big | t\right) \simeq \HGM\left(\left(\frac{1}{q},-\frac{1}{q}\right),\left(\frac{1}{r},-\frac{1}{r}\right)\Big|t^{-1}\right) ,
	\]
	and if $v_q(t-1)>0$ then $v_q(t-1) = v_q(t^{-1}-1)$, so that the roles of $q$ and $r$ are symmetric in Theorem \ref{thm:B}.
\end{remark}

\newpage
 \section{Proof of the main theorem}
 \label{Sec:proof}
 Let us start by discussing a subtle point that was overlooked in the
 literature. As explained in the introduction, Step I consists on,
 given a solution $(a,b,c)$ to (\ref{eq:GFE}), constructing a
 Galois representation attached to a rank-$2$  motive inside an abelian variety $A$.
For example, in the pioneering work of Darmon
 \cite{Darmon}, the variety $A$ is the Jacobian of the
 hyperelliptic curves studied in \cite{TTV}.

 Let $g$ denote the dimension of this variety $A$. Then, there is a totally real field $K$ of degree $g$ contained in its endomorphism algebra  
 $\End^0_K(A):=\End_K(A)\otimes \Q$. To apply
 Conjecture~\ref{conj:large-image} to $A$ one needs to determine its
 full endomorphism algebra $\End_{\overline{\QQ}}^0(A)$. If $A$ is
 absolutely simple, the condition $K \subset \End_K^0(A)$ implies that
 either $A$ has complex multiplication, or
 $\End_K^0(A) = \End_{\overline{\Q}}^0(A)$. Consequently, if $A$ does not have complex multiplication, Conjecture~\ref{conj:large-image} can be used. 

 But we remark that it is not necessarily the case that $A$ is absolutely simple. 
 In \cite[Corollary~6]{TTV} the
 authors prove 
 that for generic values, this in indeed the case. However, this does not imply that the variety
 is absolutely simple for all value of $(a,b,c)$!

 A natural solution to the problem is to replace the variety $A$
 by a \emph{building block} $B$. For this idea to work, one needs to
 have explicit control on a \emph{complete field of definition} of $B$
 (i.e.\ a field where not only the variety $B$ is defined but also all
 of its endomorphisms).  

In the rest of this section, given two abelian varieties $A_1$ and $A_2$, we will write $A_1 \sim A_2$ to denote that $A_1$ and $A_2$ are isogenous over $\overline{\Q}$.

\begin{lemma}
  \label{lemma:single-block}
  Let $A$ be an abelian variety of dimension $g$ defined over a number
  field $L$. Suppose there exists a prime ideal $\id{p}$ of $L$ where
  $A$ has toric reduction of full rank. Suppose furthermore that there
  exists a number field $K$ with $[K:\Q]=g$ and
  $K \subset \End_L^0(A)$. Then $A$ has at most one building block of
  type $I$ in Albert's classification of $\GL_2$-type, i.e.\ there
  exists a simple abelian variety $B$ defined over $\overline{\QQ}$
  and a positive integer $r$ such that $A \sim B^r$ and
  $\End_{\overline{\Q}}^0(B)$ is a totally real number field of degree $\dim(B)$.
  \end{lemma}

  \begin{proof}
    The proof follows the ideas of \cite{2412.21183}. By Poincar\'e
    decomposition theorem, there exist simple abelian varieties
    $B_1, \ldots,B_n$ and non-negative integers $r_1, \dots, r_n$ such that
    \[
      A \sim B_1^{r_1} \times \cdots \times B_n^{r_n}.
    \]
    Hence $g = r_1 \dim(B_1) + \cdots +r_n \dim(B_n)$, and
    \[
      K \subset \End_{\overline{\QQ}}^0(A) \simeq M_{r_1}(\End_{\overline{\QQ}}^0(B_1)) \times \cdots \times M_{r_n}(\End_{\overline{\QQ}}^0(B_n)).
    \]
    Since $K$ is a field, it embeds into one of the factors, say the
    first one (denote it by $B$ and $r_1$ as $r$, to ease notation). Let $H$ denote
    the center of $\End_{\overline{\QQ}}(B)$. Since $B$ is simple, this is a number field, and
    $[\End_{\overline{\QQ}}^0(B):H]=d^2$ for some positive integer
    $d$. By Albert's classification (see, e.g., \cite[Theorem 2 in \S
    21]{Mumford}), we have the following possibilities for $H$:
    \begin{itemize}
    \item type I: $H = \End_{\overline{\QQ}}^0(B)$, and so $\End_{\overline{\QQ}}^0(B)$ is a totally real
      number field and $d=1$.
    \item type II: $H$ is a totally real number field and
      $\End_{\overline{\QQ}}^0(B)$ is a totally indefinite quaternion
      algebra over $H$ (so $d=2$).
    \item type III: $H$ is a totally real number field and
      $\End_{\overline{\QQ}}^0(B)$ is a totally definite quaternion
      algebra over $H$ (so $d=2$).
    \item type IV: $H$ is a CM field, $\End_{\overline{\QQ}}^0(B)$
      is a division algebra and ${[H:\Q]}d^2 \mid 2\dim(B)$.
    \end{itemize}
    It follows from the previous description that if $M$ is any subfield
    of $\End_{\overline{\QQ}}^0(B)$ it satisfies
    \begin{equation}
    	\label{eqn:degreeineq}
      [M:\QQ] \le
      \begin{cases}
        [H:\QQ] & \text{ for type I},\\
        \dim(B) & \text{ for types II and III},\\
        2\dim(B)/d & \text{ for type IV}.
      \end{cases}
  \end{equation} 
  Since, by assumption, there exists a prime of multiplicative reduction, $B$ does not have complex multiplication (see \cite{SerreTate}) and thus does not correspond to type IV with $d=1$. 
  Since we have that
  $K \subset M_r(\End_{\overline{\QQ}}^0(B))$:
  \begin{equation}
    \label{eq:4}
     \dim(A)=g = [K:\Q] \le  r \dim(B) \le g, 
  \end{equation}
   Consequently, $A\sim B^r$.

 If $B$ is of type IV, then \eqref{eqn:degreeineq} and \eqref{eq:4} give that $d=2$ and (\ref{eq:4}) forces $K$ to be a
 maximal commutative subalgebra. Consequently, it contains the center $H$ which is a CM extension, giving a contradiction. The existence of a prime of multiplicative reduction for $B$ discards the types II and III (as 
 explained in \cite[pg.\ 126]{MR2091964}). Then
 $B$ is of type I and $[H:\QQ]=\dim(B)$ as claimed.
\end{proof}
If $A$ is an abelian variety as in the lemma, we denote by $H$ the
algebra $\End_{\overline{\Q}}^0(B)$ and by $S$ the set of prime ideals
of $L$ where $A$ has bad but non-multiplicative reduction.
\begin{prop}
  \label{prop:field-definition}
  Keep the hypotheses of Lemma \ref{lemma:single-block}. Then there exists a polyquadratic extension $M/L$
  unramified outside $S$ such that the variety $B$ and
  all of its endomorphisms are defined over $M$. 
\end{prop}

\begin{proof}
  The fact that $M$ can be taken such that $\Gal(M/L)$ is an abelian group of exponent
  $2$ follows from the work of Ribet \cite[Theorem 1.2]{MR1299737} (see also Theorem 3.24 in Guitart's Ph.D.\ thesis
  \cite{GuitartThesis}).  The ramification statement follows from
  \cite{MR371903} (see also \cite[Proposition 1.2]{MR2091964}).

\end{proof}
%
%
%
\begin{thmA}
	Let $q>r>2$ and $p>3$ be prime numbers, with $q$ and $r$ fixed.  Let $(a,b,c)$ be a hypothetical non-trivial primitive solution to
        \begin{equation*}
      	 x^q+y^p+z^r=0,          
        \end{equation*}
	and let $\rho_{\id{p}}$ be the Galois representation attached
        to $(a,b,c)$ in Theorem \ref{thm:rep-existence}. Assume that
        Conjecture \ref{conj:large-image} holds. Then there exists a
        constant $C(q,r)$ (depending only on $q$ and $r$) such that if
        $p>C(q,r)$, for any newform $f$ with complex multiplication we have
        that
        $\overline{\rho}_\id{p}\not\simeq\overline{\rho}_{f,\id{p}}$.
\end{thmA}

\begin{proof}
	Let $\bfM$ denote the maximal polyquadratic extension
	of $K$ which is unramified outside the prime ideals of $K$ dividing $qr$. Note that $\bfM/K$ is a finite extension.  Define
	\[
	C(q,r):=\max\{C(\bfM,H)\; : \; H \subset K\},
	\]
	where $C(\bfM,H)$ is the constant of Conjecture \ref{conj:large-image}.

	\vspace{5pt}
	
  Let $(a,b,c)$ be a non-trivial primitive solution to (\ref{eq:GFE}),
  where $q>r>2$ and $p>C(q,r)$, and let $K=\Q(\zeta_q)^+\cdot\Q(\zeta_r)^+$.  Let
  $\bfC$ denote the Euler curve in~(\ref{eqn:Euler}) specialized at
  $\z-spec=-\frac{a^q}{c^r}$. The curve $\bfC$ has two involutions
  (see \cite[Proposition 7.5]{GPR}), both defined over $K$. Let $\bfC'$ denote the quotient of $\bfC$ by the subgroup generated by these involutions and let $A$ be the new part of its Jacobian (so $A$ is an abelian
  variety of dimension ${(q-1)(r-1)}/{4}$). The field $K$ is
  contained in $\End_K^0(A)$ and the motive $\HGM$ of 
  (\ref{eq:motive-qr}) has a realization in $A$.

  By Catalan's conjecture \cite{MR2076124}, $b\neq \pm 1$, so there
  exists a prime number $\ell \mid b$. Let $\lambda$ be a prime ideal
  of $K$ above $\ell$. If $\lambda\nmid qr$, then the motive has
  multiplicative reduction at $\lambda$ by \cite[Theorem 1.4]{GP}. If
  $\lambda \mid qr$, since $t-1=\frac{b^p}{c^r}$, the same holds by
  Theorem~\ref{thm:B} and Remark \ref{rem:simmetry}. Then the variety
  $A$ has a prime ideal of multiplicative
  reduction.

  By \cite[Theorem 1.4]{GP}, the variety $A$ is semistable at all
  primes except possibly those dividing $qr$.  By Lemma~\ref{lemma:single-block} and
  Proposition~\ref{prop:field-definition}, there exists a simple
  $\GL_2$-type abelian variety $B$ of type I totally defined over $\bfM$
   such that $A \sim B^d$ for some $d$ (if $A$ is simple,
  $\bfM=K$ and $B=A$). Let $H=\End_{\bfM}^0(B) \subset K$.

  Let $\id{p}$ be a prime ideal of $K$ above $p$. As recalled in
  $\S \ref{Sec:reduction}$, the Galois representation $\rho_\id{p}$
  considered in Theorem \ref{thm:rep-existence} is the one attached to
  $\HGM$. Let $\id{P} = \id{p} \cap H$. Then, by definition, 
   \begin{equation}\label{eq:iso}
   \overline{\rho}_{B, \id{P}} \simeq \overline{\rho}_\id{p}|_ {\Gal(\overline{\Q}/{\bfM})}.  
   \end{equation}

    By Lemma~\ref{lemma:single-block}, the variety $B$ satisfies the hypotheses of
    Conjecture~\ref{conj:large-image}. Then, by \eqref{eq:iso}, the image of
    $\overline{\rho}_{\id{p}}$  
    restricted to $\Gal(\overline{\Q}/{\bfM})$
    contains
    $\mathrm{SL}_2(\F_{\mathfrak{p}})$ and, in particular, the same
    holds for $\overline{\rho}_{\id{p}}$.  Then,
    $\overline{\rho}_{\id{p}}$ cannot be isomorphic to the residual
    representation of any form $f$ with complex multiplication, since
    the latter has image contained in the normalizer of a Cartan
    subgroup.
	\end{proof}

\bibliographystyle{alpha}
\bibliography{biblio}
\end{document}